\title{Primes in a prescribed arithmetic progression dividing 
the sequence $\{a^k+b^k\}_{k=1}^{\infty}$}
\author{P. Moree and B. Sury}
\documentclass[12pt]{article}
\usepackage{amssymb, latexsym, amsfonts}
\def\@ptsize{2}
\newtheorem{Thm}{Theorem}
\newtheorem{Con}{Conjecture}
\newtheorem{Lem}{Lemma}
\newtheorem{cor}{Corollary}

\newtheorem{Prop}{Proposition}
\newcommand{\qed}{\hfill $\Box$}

\begin{document}
\date{}
\maketitle
%{\def\thefootnote{}
%\footnote{\noindent Max-Planck-Institut f\"ur Mathematik,
%Vivatsgasse 7, D-53111 Bonn, Germany, E-mail:
%moree@mpim-bonn.mpg.de}}
%{\def\thefootnote{}
%\footnote{{\it Mathematics Subject Classification (2001)}.
%11N37, 11R45}
\begin{abstract}
\noindent Given positive integers $a,b,c$ and $d$ such that $c$ and $d$ 
are coprime we show that the primes $p\equiv c({\rm mod~}d)$
dividing $a^k+b^k$ for some $k\ge 1$ have a natural density and explicitly
compute this density. We demonstrate our results by considering some
claims of Fermat that he made in a 1641 letter to Mersenne.
\end{abstract}
\section{Introduction}
If $S$ is a sequence of integers, then we say that an integer $m$ divides the
sequence if it divides at least one term of the sequence. The sequence
$\{a^k+b^k\}_{k=1}^{\infty}$ we will denote by $S_{a,b}$.
Several authors studied the problem of characterising
(prime) divisors of the sequence $S_{a,b}$. Hasse \cite{Hasse} seems to have
been the first to consider the Dirichlet density of prime divisors of such
sequences. Later authors, e.g., Odoni \cite{Odoni}
and Wiertelak strengthened the analytic aspects of his work, with
the strongest result being due to Wiertelak \cite{Wiertelak}. In particular, Theorem 2 of
Wiertelak \cite{Wiertelak}, in the formulation of \cite{Morseq}, yields the following corollary  
(recall that Li$(x)=\int_2^x dt/\log t$ denotes the logarithmic integral):
\begin{Thm}
\label{oud}
Let $a$ and $b$ be positive integers with $a\ne b$. 
Let $N_{a,b}(x)$ count the number of primes $p\le x$ that divide $S_{a,b}$.
Put $r=a/b$. Let $\lambda$ be the largest integer such that $r=u^{2^{\lambda}}$, with
$u$ a rational number. Let $L=\mathbb Q(\sqrt{u})$. We have
$$N_{a,b}(x)=\delta(r){\rm Li}(x)+O\left({x(\log \log x)^4\over \log^3 x}\right),$$
where the implied constant may depend on $a$ and $b$, and $\delta(r)$ is a
positive rational number that is given
in Table {\rm 0}.
\end{Thm}
\medskip
\centerline{{\bf Table 0:} The value of $\delta(r)$}
\medskip
\medskip
\begin{center}
\begin{tabular}{|c|c|c|}
\hline
$L$  & $\lambda$ & $\delta(r)$ \\
\hline
$L\ne \mathbb Q(\sqrt{2})$  & $\lambda\ge 0$ & $2^{1-\lambda}/3$ \\
\hline
$L=\mathbb Q(\sqrt{2})$  & $\lambda=0$ & $17/24$ \\
\hline
$L=\mathbb Q(\sqrt{2})$  & $\lambda=1$ & $5/12$ \\
\hline
$L=\mathbb Q(\sqrt{2})$  & $\lambda\ge 2$ & $2^{-\lambda}/3$ \\
\hline
\end{tabular}
\end{center}
\medskip
Theorem \ref{oud} implies that if $a$
and $b$ are positive integers such 
 that $a\ne b$, then asymptotically $N_{a,b}(x)\sim \delta(r)x/\log x$ with
$\delta(r)>0$. In particular, 
the set of prime divisors of the sequence $\{a^k+b^k\}_{k=1}^{\infty}$ has a positive
natural density.

In this paper we will establish, inspired by a letter from Fermat (see
next section), a related result.
\begin{Thm}
\label{main}
Let $a,b,c,d$ be positive integers with $(c,d)=1$
and assume that $a \ne b$. Let $r$ and $\lambda$ be as in the previous theorem.
Let 
$$N_{a,b}(c,d)(x):=\#\{p\le x:~p|S_{a,b},~p\equiv c({\rm mod~}d)\}.$$
Then, for 
$$ab\le \log^{2/3} x{~and~}d\le {\log^{1/6}x\over \log \log x},$$
we have
$$N_{a,b}(c,d)(x)=\delta_{a,b}(c,d){\rm Li}(x)+ O\left({2^{\lambda}x\log \log 
x\over \log^{7/6} x}\right),$$
where $\delta_{a,b}(c,d)$
is a rational number that is given
in Tables {\rm 1} to {\rm 6} and the implied constant is absolute. 
\end{Thm}
 We have
$0\le \delta_{a,b}(c,d)\le 1/\varphi(d)$ by the prime number theorem for 
arithmetic progressions. In case $\delta_{a,b}(c,d)=0$ 
there could potentially be infinitely many primes 
$p\equiv c({\rm mod~}d)$
dividing $S_{a,b}$. However, using
elementary arguments not going beyond quadratic reciprocity, one
can show that there are at most finitely many primes $p$ dividing
$S_{a,b}$ in this case. Likewise if
$\delta_{a,b}=1/\varphi(d)$, using
elementary arguments not going beyond quadratic reciprocity, one
can show that in each case there are at most finitely many primes 
$p\equiv c({\rm mod~}d)$
not 
dividing $S_{a,b}$.
For a more precise statement we refer to Theorem \ref{moetnaamhebben}.

Inspection of the tables shows that we can always write
$\varphi(d)\delta_{a,b}(c,d)={c\over 2^m\cdot 3}$, for some 
non-negative integers $c$ and $m$.\\

\noindent
{\bf Notations:}\\
As the tables for the density depend on 
some auxiliary parameters computed from $a,b,c,d$, 
some notations are needed to read them. We introduce these notations here and
they will be maintained throughout this article. 
Given $a,b$ and the modulus $d$, there is a unique table among the 6 from 
which one reads off the density.  
Put $r = a/b = r_0^h$, where $r_0$ is not a proper power of a rational 
number. Write $h=2^{\lambda}h', d = 2^{\delta}d'$, with $h',d'$ odd. Put 
$v_2(c-1)=\gamma$, where it is understood that $\gamma$ is larger than any 
number when $c=1$. We denote the discriminant of the quadratic field 
$\mathbb{Q}(\sqrt{t})$ by $D(t)$ and we put $D(r_0) = 
2^{\delta_0}D'$. We also write $r_0=u/v$ and $t = -r_0$ or $\prod_{i=1}^k 
(\frac{-1}{p_i})p_i$ according as to whether $uv$ is odd or 
$uv = 2 \prod_{i=1}^k p_i$. By $d^{\infty}$ we denote the supernatural
(Steinitz) number $\prod_{p|d}p^{\infty}$.
For each positive integer $j \geq 1$, we put
$N_j = \mathbb{Q}(\zeta_{2^j},r^{1/2^{j-1}},\zeta_d)$ 
and $N'_j = \mathbb{Q}(\zeta_{2^j},r^{1/2^{j}},\zeta_d)$, 
where $\zeta_l$ for any $l$, denotes any fixed primitive $l$-th root of 
unity. Finally, for $j \geq 1$, the intersection fields 
$K_j := \mathbb{Q}(\zeta_{2^j},r^{1/2^{j-1}}) \cap \mathbb{Q}(\zeta_d)$ 
and 
$K'_j := \mathbb{Q}(\zeta_{2^j},r^{1/2^{j}}) \cap \mathbb{Q}(\zeta_d)$ 
will occur throughout our discussion.
\newpage

{\bf Table 1 : $\mathbb{Q}(\sqrt{r_0}) \neq \mathbb{Q}(\sqrt{2}), D'\nmid d'$}
\vskip 3mm

\begin{tabular}{|r|c|r|}
\hline
$\lambda$ & $\delta$ & $\phi(d) \delta_{a,b}(c,d)$ \\ \hline
$< \delta$ & $\leq \gamma$ & $1 - 
\frac{2^{\lambda+1-\delta}}{3}$\\ \hline
$\ast$ & $>0, \leq {\rm min}(\lambda,\gamma)$ & $\frac{2^{\delta-\lambda}}{3}$\\ 
$\ast$ & $0$ & $\frac{2^{1-\lambda}}{3}$\\ 
\hline
$\geq \gamma$ & $> \gamma$ & $0$\\ \hline
$< \gamma$ & $> \gamma$ & 
$1 - 2^{\lambda-\gamma}$\\ \hline

\end{tabular}
\vskip 8mm

{\bf Table 2 : $\mathbb{Q}(\sqrt{r_0}) \neq \mathbb{Q}(\sqrt{2}), 
D'|d', \delta_0 \leq \delta$}
\vskip 3mm

\begin{tabular}{|r|c|r|r|}
\hline
$\lambda$ & $\delta$ & $(\frac{D(r_0)}{c})$ & 
 $\phi(d) \delta_{a,b}(c,d)$\\ 
\hline
$\geq \delta-1$ & $>0, \leq \gamma$ & $1$ & 
$\frac{2^{\delta-1-\lambda}}{3}$\\ 
 &  & $-1$ & $2^{\delta-1-\lambda}$\\ \hline
$\ast$ & $0$ & $1$ & $\frac{2^{-\lambda}}{3}$\\ 
 &  & $-1$ & $2^{-\lambda}$\\ 
\hline
$< \delta-1$ & $\leq \gamma$ & $1$ &
$1 - \frac{2^{\lambda+2-\delta}}{3}$\\ 
 & & $-1$ & $1$\\ \hline
$\geq \delta$ & $>\gamma$ & $\ast$ & $0$\\ \hline
$\leq \gamma-1$ & $> \gamma$ & $1$ & 
$1 - 2^{\lambda+1-\gamma}$\\ 
 & & $-1$ & $1$\\ \hline
$\geq \gamma$ & $> \lambda$ & $\ast$ & $0$\\ \hline

\end{tabular}
\vskip 8mm

{\bf Table 3 : $\mathbb{Q}(\sqrt{r_0}) \neq \mathbb{Q}(\sqrt{2}), 
D'|d'$ and $\delta_0 > \delta$}
\vskip 3mm

\begin{tabular}{|r|c|r|r|}
\hline
$\lambda$ & $\delta$ & $(\frac{D(t)}{c})$ & 
$\phi(d) \delta_{a,b}(c,d)$\\ 
\hline
\hline
$<\delta-1$ & $\leq \gamma$ & $1$ & 
$1- 
\frac{2^{\lambda+1-\delta}}{3} 
+ \frac{2^{\lambda+2+\delta-2\delta_0}}{3}$\\
\hline
$<\delta-1$ & $\leq \gamma$ & $-1$ & 
$1- 
\frac{2^{\lambda+1-\delta}}{3} 
- \frac{2^{\lambda+2+\delta-2\delta_0}}{3}$\\
\hline 
$ = \delta-1$ & $\leq \gamma$ & $1$ & $\frac{2}{3} 
+ \frac{2^{2\delta+1-2\delta_0}}{3}$\\
\hline
$ = \delta-1$ & $\leq \gamma$ & $-1$ & $\frac{2}{3} 
-\frac{2^{2\delta+1-2\delta_0}}{3}$\\
\hline
$\leq \gamma-1$ & $>\gamma$ & $\ast$ & 
$1-2^{\lambda-\gamma}$\\ \hline
$\geq \gamma$ & $> \lambda$ & $\ast$ & $0$\\ \hline
$\geq \delta$ & $>\gamma$ & $\ast$ & $0$\\ \hline
$\leq \delta_0-2$ & $>0,\leq {\rm min}(\gamma,\lambda)$ & $1$ & 
$\frac{2^{\delta-\lambda}}{3} + 
\frac{2^{\lambda+2+\delta-2\delta_0}}{3}$\\ \hline
$\leq \delta_0-2$ & $>0, \leq {\rm min}(\gamma,\lambda)$ & $-1$ & 
$\frac{2^{\delta-\lambda}}{3} - 
\frac{2^{\lambda+2+\delta-2\delta_0}}{3}$\\
\hline
$\geq \delta_0-1$ & $>0,\leq \gamma$ & $1$ &
$\frac{2^{\delta-1-\lambda}}{3}$\\ \hline
$\geq \delta_0-1$ & $>0,\leq \gamma$ & $-1$ &
$2^{\delta-\lambda-1}$\\
\hline
$\leq \delta_0-2$ & $0$ & $1$ & 
$\frac{2^{1-\lambda}}{3} + 
\frac{2^{\lambda+3-2\delta_0}}{3}$\\ \hline
$\leq \delta_0-2$ & $0$ & $-1$ & 
$\frac{2^{1-\lambda}}{3} - 
\frac{2^{\lambda+3-2\delta_0}}{3}$\\
\hline
$\geq \delta_0-1$ & $0$ & $1$ &
$\frac{2^{-\lambda}}{3}$\\ \hline
$\geq \delta_0-1$ & $0$ & $-1$ &
$2^{-\lambda}$\\
\hline

\end{tabular}
\vfil\eject

{\bf Table 4 : $\mathbb{Q}(\sqrt{r_0}) = \mathbb{Q}(\sqrt{2}), 
\delta \leq 2$}
\vskip 5mm

\begin{tabular}{|r|r|r|r|}
\hline
$\lambda$ & $\delta$ & $\gamma$ & 
$\phi(d) \delta_{a,b}(c,d)$
\\ 
\hline \hline
$0$ & $\leq 1$ & $\geq \delta$ & $17/24$\\ \hline
$0$ & $2$ & $\geq \delta$ & $11/12$\\ \hline
$0$ & $2$ & $1$ & $1/2$\\ \hline
$1$ & $2$ & $1$ & $0$\\ \hline
$1$ & $\leq 1$ & $\geq \delta$ & $5/12$\\ \hline
$1$ & $2$ & $\geq \delta$ & $5/6$\\ \hline
$\ge 2$ & $\le 1$ & $\geq \delta$ & $2^{-\lambda}/3$\\ \hline
$\ge 2$ & $2$ & $\geq \delta$ & $2^{1-\lambda}/3$\\ \hline
$\ge 2$ & $2$ & $1$ & $0$\\ \hline

\end{tabular}
\vskip 5mm

{\bf Table 5 : $\mathbb{Q}(\sqrt{r_0}) = \mathbb{Q}(\sqrt{2}), 
\delta \geq 3, \lambda > 0$}
\vskip 5mm

\begin{tabular}{|r|r|r|r|}
\hline
$\lambda$ & $\delta$ & $\gamma$ & 
$\phi(d) \delta_{a,b}(c,d)$\\
\hline \hline
$\geq 2$ & $3$ & $< \delta$ & $0$\\ \hline
$\geq \delta-1$ & $\geq 3$ & $\geq \delta$ & 
$\frac{2^{\delta-1-\lambda}}{3}$\\ 
\hline
$\geq 2, < \delta-1$ & $\geq 4$ & $\geq \delta$ & $1 - 
\frac{2^{\lambda+2-\delta}}{3}$\\ \hline
$\geq 2, \leq \gamma-2$ & $\geq 4$ & $< \delta$ & 
$1-2^{\lambda+1-\gamma}$\\ \hline
$\geq {\rm max}(2,\gamma-1)$ & $\geq 4$ & $< \delta$ & $0$\\ \hline
$1$ & $\geq 3$ & $\geq \delta$ & $1- \frac{2^{3-\delta}}{3}$\\ 
\hline
$1$ & $\geq 3$ & $1$ & $0$\\ \hline
$1$ & $\geq 3$ & $2$ & $1$\\ \hline
$1$ & $\geq 3$ & $>3,<\delta$ & $1-2^{2-\gamma}$\\ 
\hline \hline

\end{tabular}
\vskip 10mm

{\bf Table 6 : $\mathbb{Q}(\sqrt{r_0}) = \mathbb{Q}(\sqrt{2}), 
\delta \geq 3, \lambda = 0$}
\vskip 5mm

\begin{tabular}{|r|r|r|}
\hline
 $\gamma$ & $c({\rm mod~}8)$ & $\phi(d) \delta_{a,b}(c,d)$\\ 
\hline \hline
 $\geq \delta$ & $1$ & 
$1-\frac{2^{2-\delta}}{3}$\\ \hline
 $\leq 2$ & $\pm{1}$ & 
$0$\\ \hline
 $\leq 2$ & $\pm{3}$ & 
$1$\\ \hline
 $\geq 3,< \delta$ & $1$ & 
$1-2^{1-\gamma}$\\ \hline

\end{tabular}
\vskip 5mm
\noindent In the next section we reconsider a letter from Fermat and papers by 3 authors
\cite{aigner, brauer, sier} in the light of Theorem \ref{main}. In Section 3 we  prove
Theorem \ref{main}, except for the fact that an expression
for $\delta_{a,b}(c,d)$ in terms of data from algebraic number theory appears. 
In Sections 4-7 we evaluate
this expression for $\delta_{a,b}(c,d)$. The outcome is recorded in Tables 1-6. 
This then completes the proof of Theorem \ref{main}.
In
Section 8 we determine the cases in which $\delta_{a,b}(c,d)=0$, respectively 
$\delta_{a,b}(c,d)=1/\varphi(d)$.
In the final section we produce the results of some numerical experiments and show 
that they match well with what can be read from our tables.

\section{On a letter of Fermat to Mersenne}
Fermat \cite[p. 220]{Fermat}, cf. 
Dickson \cite[p. 267]{Dickson}, in a letter to Mersenne dated 15 June 1641 stated that ($p$ will always be used to denote primes): 
\begin{Con}{\rm (Fermat, 1641)}\\
{\rm 1)} If $p|S_{3,1}$, then $p\not\equiv -1({\rm mod~}12)$.\\
{\rm 2)} If $p|S_{3,1}$, then $p\not\equiv +1({\rm mod~}12)$.\\
{\rm 3)} If $p|S_{5,1}$, then $p\not\equiv -1({\rm mod~}10)$.\\
{\rm 4)} If $p|S_{5,1}$, then $p\not\equiv  +1({\rm mod~}10)$.
\end{Con}
Pur $r=a/b$. For $p\nmid ab$ there exists a smallest positive integer $k$ such
that $r^k\equiv 1({\rm mod~}p)$; this is ord$_p(r)$, the multiplicative order
of $r({\rm mod~}p)$. It is not difficult to see that if $p\nmid ab$, then $p|S_{a,b}$
if and only if
ord$_p(r)$ is even. If $p|ab$ and $p\nmid (a,b)$, then clearly $p\nmid S_{a,b}$. (With $(a,b)$
and $[a,b]$ we denote the greatest common divisor, respectively lowest common multiple of
$a$ and $b$.)
Using this
observation and the law of quadratic reciprocity it is easy to see that
the following holds:
\begin{Prop}
\label{propje}
Conjecture 1.1 of Fermat holds true.
\end{Prop}
{\it Proof}. For $p>3$ by the law
of quadratic reciprocity we have $({3\over p})({p\over 3})=(-1)^{p-1\over 2}$.
Suppose that $p\equiv -1({\rm mod~}12)$. It then follows that $({3\over p})=1$.
By Euler's identity we then have
$3^{p-1\over 2} \equiv ({3\over p})=1({\rm mod~}p)$. Since $(p-1)/2$ is the
largest odd divisor of $p-1$ it follows that ord$_p(3)$ is odd. This implies
that $p\nmid S_{3,1}$. \qed\\

\noindent However, a computeralgebra computation learns that the remaining conjectures
are all false. Counterexamples (in ascending order) are listed below :\\
Counterexamples to:\\

\noindent Conjecture 1.2: 37, 61, 73, 97, 157, 193, 241, 337, 349, 373, 397, 409, 457, $\cdots $\\
Conjecture 1.3: 41, 61, 241, 281, 421, 521, 601, 641, 661, 701, 761, 821, 881,$\cdots $\\
Conjecture 1.4: 29, 89, 229, 349, 449, 509, 709, 769, 809, 929, 1009, 1049, $\cdots $\\

\noindent Sierpi\'nski 
suggested that Conjecture 1.2 is false for infinitely many primes. This was
proved by Schinzel \cite{schinzel}, who in the same paper showed that also
Conjecture 1.3 and Conjecture 1.4 are false for infinitely many primes.
Theorem \ref{main} implies that there is even a positive density
of primes for which the conclusions of these three conjectures are false:
\begin{cor}
We have $$\delta_{3,1}(1,12)={1\over 6},~\delta_{3,1}(5,12)={1\over 4},~\delta_{3,1}(7,12)=
{1\over 4}{\rm ~and~}\delta_{3,1}(11,12)=0.$$
Furthermore, we have
$$\delta_{5,1}(1,10)={1\over 12},~\delta_{5,1}(3,10)={1\over 4},~\delta_{5,1}(7,10)=
{1\over 4}{\rm ~and~}\delta_{5,1}(9,10)={1\over 12}.$$
In particular, the relative density of the primes for which the conclusion in
Conjectures 1.1-1.4 fail are, respectively,
$${\delta_{3,1}(11,12)\over \delta(3)}=0,~
{\delta_{3,1}(1,12)\over \delta(3)}={1\over 4},~
{\delta_{5,1}(9,10)\over \delta(5)}={1\over 8},~{\delta_{5,1}(1,10)\over \delta(5)}={1\over 8}.$$ 
\end{cor}
After Fermat various authors considered primes in arithmetic progressions
dividing $S_{a,b}$. Thus  Sierpi\'nski \cite{sier} proved that  
every prime $p\equiv \pm 3({\rm mod~}8)$ divides $S_{2,1}$ and, furthermore,
that no prime $p\equiv 7({\rm mod~}8)$ divides $S_{2,1}$. This result
easily follows on using that $({2\over p})=(-1)^{(p^2-1)/8}$. Sierpi\'nski states
that M.A. Makowski has proved that infinitely many primes 
$p\equiv 1({\rm mod~}8)$ divide $S_{2,1}$ (namely Makowski notices that
the prime factors of the numbers of the form $2^{2^n}+1$ with $n\ge 3$ have
the required property) and ends his paper with stating the problem of whether there
are infinitely many primes $p\equiv 1({\rm mod~}8)$ not dividing $S_{2,1}$.
Subsequently, using results on the biquadratic and octavic residue
character of 2, this problem has been independently resolved by A. Aigner \cite{aigner}
and A. Brauer \cite{brauer}. Brauer shows for example that the infinitely
many primes $p\equiv 9({\rm mod~}16)$ which can be represented as
$65x^2+256xy+256y^2$ all do not divide $S_{2,1}$ (the number
of such primes $\le x$ is of order $O(x/\sqrt{\log x})$ by a result
of G. Pall \cite{pall}, and thus
this set has natural density zero). Using the first entry
of Table 6 we infer that there many more primes
not dividing $S_{2,1}$: $1/6$th of all primes $p\equiv 1({\rm mod~}8)$
do not divide $S_{2,1}$.

\section{The density written as infinite sum}
In order to evaluate $\delta_{a,b}(c,d)$ we will make use of the following result.
\begin{Thm}
\label{een}
Let $a,b,c,d$ be positive integers with $c\ge 1$ and $d\ge 1$ coprime. Let $\sigma_c$ denote
the automorphism of $\mathbb Q(\zeta_d)$ determined by $\sigma_c(\zeta_d)=
\zeta_d^c$. 
The density $\delta_{a,b}(c,d)$ of primes $p\equiv c({\rm mod~}d)$ such that
$p|S_{a,b}$  exists
and satisfies
\begin{equation}
\label{begin}
\delta_{a,b}(c,d)=\sum_{j=1}^{\infty}
\Big({\tau(j)\over [N_j:\mathbb Q]}
-{\tau'(j)\over [N'_j:\mathbb Q]}\Big),
\end{equation}
where $$\tau(j)=\cases{1 & if $\sigma_c|_{K_j}=$id.;\cr 0 & otherwise,}
{\rm ~and,~similarly,~}\tau'(j)=\cases{1 & if $\sigma_c|_{K'_j}=$id.;\cr 0 & otherwise.}$$
Furthermore, Theorem \ref{main} holds true with $\delta_{a,b}(c,d)$ as given by (\ref{begin}).
\end{Thm}
{\it Proof}. In case ord$_p(r)$ is defined we can define the 
{\it index},
$i_p(r)$, as $(p-1)/{\rm ord}_p(r)$. Note that it equals 
$[\mathbb F_p^*:\langle r\rangle]$.
There is a unique $j\ge 1$ such that $2^{j-1}||i_p(r)$.
Let $P_j$ denote the set of primes $p$ such that $2^{j-1}||i_p(r)$.
Note that $\cup_{j=1}^{\infty}P_j$ equals, with finitely many exceptions, 
the set
of all primes and that the $P_i$ are disjoint sets. Now note that for a 
prime $p$ in $P_j$
we have that ord$_p(r)$ is even if and only if $p\equiv 1({\rm mod~}2^j)$. Thus, 
except for finitely many
primes, the set of prime divisors of $S_{a,b}$ satisfying 
$p\equiv c({\rm mod~}d)$ is of the form $\cup_{j=1}^{\infty}Q_j$, where
$$Q_j:=\{p:p\equiv c({\rm mod~}d),p\equiv 1({\rm mod~}2^j),~p\in P_j\}.$$
It is an easy observation that $n|i_p(r)$ if and only if $p$ splits completely in 
$\mathbb{Q}(\zeta_n,r^{1/n})$.
Using this observation and writing `s.c.' below to mean that the prime is 
split completely, we infer that 
$$Q_j=\{p:p\equiv c({\rm mod~}d),p~{\rm s.c. in~}
\mathbb{Q}(\zeta_{2^j},r^{1/2^{j-1}}),{\rm but~not~s.c. in~}
\mathbb{Q}(\zeta_{2^j},r^{1/2^{j}}) \}.$$
On invoking the Chebotarev density theorem, it is then found that the set 
$Q_j$ has a natural
density that is given by
$$\delta(Q_j)={\tau(j)\over [N_j:\mathbb Q]}
-{\tau'(j)\over [N'_j:\mathbb Q]}.$$
On proceeding as in the proof of Lemma 8 of \cite{ulmer} it is then found that
for $ab\le \log^{2/3}x$ and $[d,2^j]\le y:=\log^{1/6}x/\log \log x$, and any number
$A>0$, we have
\begin{equation}
\label{zeur1}
Q_j(x)=\delta(Q_j){\rm Li}(x)+O_A\Big({x\over \log^A x}\Big).
\end{equation}
Thus
$$N_{a,b}(c,d)(x)=\sum_{j\ge 1}Q_j(x)=\sum_{[d,2^j]\le
y}Q_j(x)+O(\sum_{[d,2^j]>y}\pi(x;[2^j,d],c_j)),$$
where $\pi(x;m,n)$ denotes the number of primes $p\le x$ such that $p\equiv n({\rm
mod~}m)$ and $c_j$ is any integer such that $c_j\equiv c({\rm mod~}d)$ and $c_j\equiv 1({\rm mod~}2^j)$
if such an integer exists and 1 otherwise. A minor modification of the proof of
Lemma 2 of \cite{polen} then yields that 
\begin{equation}
\label{zeur2} 
N_{a,b}(c,d)(x)=\sum_{[d,2^j]\le
y}Q_j(x)+O\Big({x\log \log x\over \log^{7/6}x}\Big).
\end{equation}
Using Lemma \ref{tweegraad} we find that
\begin{equation}
\label{zeur3}
\sum_{[d,2^j]>y}^{\infty}\delta(Q_j)=O\Big(2^{\lambda}\sum_{[d,2^j]>y}{1\over [d,2^j]2^j}\Big)
=O({2^{\lambda}\over y}).
\end{equation}
On combining (\ref{zeur1}), (\ref{zeur2}) and (\ref{zeur3}),
the result is then obtained with
$\delta_{a,b}(c,d)=\sum_{j=1}^{\infty}\delta(Q_j)$.\qed\\

\noindent 
{\tt Remark 1}. The algebraic side of the approach above 
(originating in Moree \cite{polen})
is not the traditional one to study the divisiblity
of sequences $S_{a,b}$, but is chosen since it turns out to be easier to explicitly work
out. The traditional approach rests on the observation that
if $p\equiv 1+2^j({\rm mod~}2^{j+1})$ for some $j$ (which is uniquely determined), then
ord$_p(r)$ is odd if and only if $r^{(p-1)/2^j}\equiv 1({\rm mod~}p)$, that is if and
only if $p$ splits completely
in $\mathbb{Q}(\zeta_{2^j},r^{1/2^{j}})$, 
see e.g. \cite{morsurvey} for a sketch of the traditional approach. Note
that $(p-1)/2^j$ is the largest odd divisor of $p-1$ and so ord$_p(r)$ is odd if and only
if
ord$_p(r)$ divides $(p-1)/2^j$.\\

\noindent {\tt Remark 2}. On GRH the existence of $\delta_{a,b}(c,d)$ was established by Moree \cite[Theorem 1]{Mor2}. He 
showed under GRH that the set of primes $p$ such that $p\equiv a_1({\rm mod~}d_1)$ and
ord$_p(r)\equiv a_2({\rm mod~}d_2)$ has a density $\delta_r(a_1,d_1;a_2;d_2)$ and gave an 
expression for it in terms
of field degrees and Galois intersection coefficients ($\tau(j)$ and $\tau'(j)$ in 
Theorem \ref{een} are
examples of such coefficients). Since $\delta_{a,b}(c,d)=\delta_r(c,d;0,2)$, 
where $r=a/b$, it follows
that  $\delta_{a,b}(c,d)$ exists under GRH.\\

\noindent From our tables it is seen that $\delta_{a,b}(c,d)$ is always rational. Below
a conceptual explanation for this is given.
\begin{Prop}
The density $\delta_{a,b}(c,d)$ is always a rational number.
\end{Prop}
{\it Proof}. We show that the sum in (\ref{begin}) always yields a rational number. Note that
$K_j\subseteq K_{j+1}$ and $K_j'\subseteq K'_{j+1}$ and hence the fields
$\lim_{j\rightarrow \infty}K_j,~\lim_{j\rightarrow \infty}K'_j$ exist. Denote these
limits by $K,K'$. Note that $K=K'$. It follows that there exists $j_0$ such that
$\tau(j)=\tau'(j)$ and $K_j=K'_j=K=K'$ for every $j\ge j_0$. By Lemma \ref{tweegraad} it follows
that there exist constants $c_1$ and $c_2$ such that $[N_j:\mathbb Q]=c_14^j$ and
$[N'_j:\mathbb Q]=c_2 4^j$ for every $j$ large enough. It follows that the terms with
$j$ large enough in (\ref{begin}) are in geometric progression and sum to a rational number. 
The terms are all rational and so $\delta_{a,b}(c,d)$ is itself rational.\qed

\section{Preliminaries on field degrees and field intersections}
The following facts from elementary algebraic number theory,  
for further details we refer to e.g. Moree \cite{Mor2}, will be
used freely in the sequel: \\
1) a quadratic field $K\subseteq \mathbb Q(\zeta_n)$
iff the discriminant of $K$ divides $n$.\\
2) Let $\mathbb Q(\sqrt{\Delta})\subseteq \mathbb Q(\zeta_n)$ be a quadratic
fields of discriminant $\Delta$ and $b$ be an integer with $(b,n)=1$. Then
$\sigma_b|_{\mathbb Q(\sqrt{\Delta}}={\rm id.}$ iff $({\Delta\over b})=1$, with
$({\cdot\over \cdot})$ the Krnecker symbol.\\
\indent In order to use Theorem \ref{een} 
to compute $\delta_{a,b}(c,d)$, we first compute the 
degrees of the fields $N_j,N'_j$ for $j \geq 1$. This can be done 
directly or by using the general formula from Lemma 1 of \cite{Mor1} 
quoted below:
\begin{Lem}
Put $n_t=[2^{v_2(ht)+1},D(r_0)]$. We have
$$[\mathbb{Q}(\zeta_{kt},r^{1/k}):\mathbb{Q}]={\phi(kt)k\over 
\epsilon(kt,k)(k,h)},~{\rm ~where~}\epsilon(kt,k)=\cases{2 & if $n_t|kt$;\cr 1 & if $n_t\nmid kt$.}$$
\end{Lem}
Using the lemma or otherwise, we compute the degrees of 
$$\cases{N_j = \mathbb{Q}(\zeta_{2^j},r^{1/2^{j-1}},\zeta_d) =
\mathbb{Q}(\zeta_{2^{{\rm max}(j,\delta)}d'},r^{1/2^{j-1}}); & \cr
N'_j = 
\mathbb{Q}(\zeta_{2^j},r^{1/2^j},\zeta_d) 
= \mathbb{Q}(\zeta_{2^{{\rm max}(j,\delta)}d'},r^{1/2^j}), & \cr}$$
to be as given in Lemma \ref{tweegraad}.
The degrees turn out to be dependent on the following property which we 
call $C_j$ :\\

\noindent {\it The property $(C_j)$ holds if and only if} 
$D'|d' , \delta_0 \leq {\rm max}(j, \delta).$\\

\noindent Note that if $D'|d'$, then $(C_j)$ can fail only for finitely many $j$'s.
\vskip 3mm

\begin{Lem}
\label{tweegraad}
The degrees of $N_j = \mathbb{Q}(\zeta_{2^j},r^{1/2^{j-1}},\zeta_d)$
and $N'_j = \mathbb{Q}(\zeta_{2^j},r^{1/2^{j}},\zeta_d)$ over $\mathbb{Q}$ 
are given by:
$${1\over \varphi(d)}[N_j: \mathbb{Q}]=\cases{ 
 2^{{\rm max}(j,\delta)-1}  & if $j\leq \lambda + 1$;\cr
 2^{{\rm max}(j,\delta)+j-\lambda-3} & if
$j > \lambda+1$ and $(C_j)$ holds;\cr
2^{{\rm max}(j,\delta)+j-\lambda-2}  & if
$j> \lambda+1,$ and $(C_j)$ fails,\cr}$$
$${1\over \varphi(d')}[N'_j: \mathbb{Q}]=\cases{2^{{\rm max}(j,\delta)-1}  & if $j \leq \lambda$;\cr
2^{{\rm max}(j,\delta)+j-\lambda-2}  & if $j> \lambda$ and $(C_j)$ holds;\cr
2^{{\rm max}(j,\delta)+j-\lambda-1}  & if $j>\lambda$ and $(C_j)$ fails.}$$
\end{Lem}
\vskip 3mm

\noindent
{\tt Remark 3}. {\it Equivalent form of $(C_j)$}.\\
It will also be convenient to use the following version of $(C_j)$ later. 
\vskip 2mm

\noindent
Property {\it $(C_j)$ holds if and only if,
either $D(r_0)|d$ or $D(r_0)|2^ld, D(r_0) \nmid 2^{l-1}d$ for some $l \geq 
1$ and $j \geq l+\delta$.\\
Equivalently, property $(C_j)$ fails if, and only if,
either $D(r_0) \nmid 2^ld~ \forall l \geq 0$ or
$D(r_0)|2^ld, D(r_0) \nmid 2^{l-1}d$ for some $l \geq 1$ and
$j< l+ \delta.$}
\vskip 5mm

\noindent
In the remainder of this section we assume that $\mathbb{Q}(\sqrt{r_0}) \neq 
\mathbb{Q}(\sqrt{2})$. The case $\mathbb{Q}(\sqrt{r_0}) \neq 
\mathbb{Q}(\sqrt{2})$
requires modification due to the ramification of $2$ in 
cyclotomic extensions generated by large $2$-power roots of unity
and is discussed in Sections 7 and 8.

We need to determine precisely the set of all $j \geq 1$
for which $\tau(j)=1$ and those for which $\tau'(j)=1$. To this end
we first determine the degrees of $K_j,K'_j$ over $\mathbb{Q}$.
\vskip 3mm

\begin{Lem}
\label{vvier}
When $\delta>0$, the degrees of $K_j,K'_j$ are given by the expressions :
\vskip 3mm

$$[K_j : \mathbb{Q}]=\cases{2^{{\rm min}(j,\delta)} & if $j \leq \lambda+1$;\cr
2^{{\rm min}(j,\delta)} & if $j > \lambda+1$ and $(C_j)$ holds;\cr
2^{{\rm min}(j,\delta)-1} & if $j > \lambda+1$ and $(C_j)$ does not hold,}$$

$$[K'_j : \mathbb{Q}]=\cases{
2^{{\rm min}(j,\delta)} & if $j \leq \lambda$;\cr
2^{{\rm min}(j,\delta)} & if $j > \lambda$ and $(C_j)$ holds;\cr
2^{{\rm min}(j,\delta)-1} & if $j>\lambda$ and $(C_j)$ does not hold.}$$

\end{Lem}
\vskip 3mm

\noindent
{\it Proof}. When $j\leq \lambda + 1$, clearly  $r^{1/2^{j-1}}$ is 
rational
and, therefore, $K_j = \mathbb{Q}(\zeta_{2^{{\rm min}(j,\delta)}})$.  
Similarly,
$K'_j = \mathbb{Q}(\zeta_{2^{{\rm min}(j,\delta)}})$ if $j \leq 
\lambda$. Further, note that $K_j \subseteq K'_j$ for all $j$. Writing 
$L_j = \mathbb{Q}(\zeta_{2^j},r^{1/2^{j-1}})$, and
$L'_j = \mathbb{Q}(\zeta_{2^j},r^{1/2^j})$, we have
$N_j = L_j \mathbb{Q}(\zeta_d)$ and $K_j = 
L_j \cap \mathbb{Q}(\zeta_d)$. Therefore, 
$$[K_j: \mathbb{Q}] = \frac{[L_j: \mathbb{Q}][\mathbb{Q}(\zeta_d): 
\mathbb{Q}]}{[N_j: \mathbb{Q}]}.$$
Similarly, $N'_j = L'_j \mathbb{Q}(\zeta_d)$ and $K'_j = 
L'_j \cap \mathbb{Q}(\zeta_d)$. So, 
$$[K'_j: \mathbb{Q}] = \frac{[L'_j: \mathbb{Q}][\mathbb{Q}(\zeta_d): 
\mathbb{Q}]}{[N'_j: \mathbb{Q}]}.$$
Using the above degree computations for $N_j,N'_j$ etc., we obtain
the asserted expressions.
\vskip 5mm

\noindent
For $\delta=0$, the above formula has to be modified as we have used 
$\phi(2^{\delta}) = 2^{\delta-1}$. In this case, we get :
\vskip 3mm

\begin{Lem}
When $\delta= 0$, we have 
$$[K_j : \mathbb{Q}]=\cases{2 & if $j > \lambda+1$ and $(C_j)$ holds;\cr
1 & if either $j \leq \lambda+1$ or $j > \lambda+1$ and $(C_j)$ fails,}$$
and
$$[K'_j : \mathbb{Q}]= \cases{2 & if $j > \lambda$ and $(C_j)$ holds;\cr
1 & if either $j \leq \lambda$ or $j>\lambda$ and $(C_j)$ fails.}$$
\end{Lem}

\noindent
{\tt Remark 4}. Since $K_j$ is a subfield of $K'_j$, it follows from the above degree 
computation that
$K_j = K'_j$ in all cases except possibly when $j = \lambda +1$.
For $j =\lambda+1$, we have 
$\mathbb{Q}(\zeta_{2^{{\rm min}(\lambda+1,\delta)}}) = K_{\lambda+1}$ 
and the 
degree of $K'_{\lambda+1}$ 
over $K_{\lambda+1}$ is $2$ if $D'|d'$ and $\delta_0 \leq {\rm max}(\lambda+1, 
\delta)$. If 
this latter condition $(C_{\lambda+1})$ does not hold, 
then $K_{\lambda+1} = K'_{\lambda+1}$. In other words, we have the following property :
$$K_j=K'_j,~~~\tau(j)=\tau'(j)~~~ \forall~~~j \neq \lambda+1.$$
We would like to actually write the fields $K_j,K'_j$ in a convenient 
form so that we can determine how the automorphism $\zeta_d 
\mapsto \zeta_d^c$ acts on them. 
Note that clearly the field 
$\mathbb{Q}(\zeta_{2^{{\rm min}(j,\delta)}})$
is always contained in $K_j,K'_j$ and its degree 
is either the whole or half of that of $K_j, K'_j$ .
We look for a subfield of the form 
$\mathbb{Q}(\zeta_{2^{{\rm min}(j,\delta)}})$
or 
$\mathbb{Q}(\zeta_{2^{{\rm min}(j,\delta)}}, \sqrt{v})$ which has the 
full degree and will, therefore, have to be the whole field. 
\vskip 5mm

\begin{Lem}
For $j \leq \lambda$, $K_j = K'_j = 
\mathbb{Q}(\zeta_{2^{{\rm min}(j,\delta)}})$.\\
Furthermore, $K_{\lambda+1} =  \mathbb{Q}(\zeta_{2^{{\rm min}(\lambda+1,\delta)}})$.\\
For $j> \lambda+1$,  $K_j = K'_j$.\\
For $j\geq \lambda+1$,  $K'_j$ is :\\
{\rm (a)} $\mathbb{Q}(\zeta_{2^{{\rm min}(j,\delta)}})$ if either $D' \nmid d'$ 
or if $\delta_0 > {\rm max}(j, \delta)$ ;\\
{\rm (b)} $\mathbb{Q}(\zeta_{2^{{\rm min}(j,\delta)}}, \sqrt{r_0})$ if 
$D(r_0)|d$ ;\\
{\rm (c)} $\mathbb{Q}(\zeta_{2^{{\rm min}(j,\delta)}}, \sqrt{-r_0})$
if $D'|d', \delta < \delta_0 \leq {\rm max}(j, \delta)$, 
where $r_0 = u/v$ and $2\nmid uv$;\\
{\rm (d)} $\mathbb{Q}(\zeta_{2^{{\rm min}(j,\delta)}}, \sqrt{\prod_{i=1}^k 
(\frac{-1}{p_i})p_i})$  if $D'|d', \delta < \delta_0 \leq {\rm max}(j, \delta)$, 
where $r_0=u/v$ with $uv = 2 \prod_{i=1}^k p_i$ and $p_i>2$ for $i=1,\ldots,k$.
\end{Lem}
{\it Proof}. We know that $K_j = K'_j = \mathbb{Q}(\zeta_{2^{{\rm 
min}(j,\delta)}})$ if either $j \leq \lambda$ or $j> \lambda+1$ and 
$(C_j)$ fails. Also,
$K_{\lambda+1} = \mathbb{Q}(\zeta_{2^{{\rm min}(\lambda+1,\delta)}}) = 
K'_{\lambda+1}$ unless $(C_{\lambda+1})$ fails. In other words, we have to 
determine $K'_j$ only for those $j > \lambda$ for which $(C_j)$ holds.\\
Recall that the truth of $(C_j)$ is equivalent to the property :\\
either $D(r_0)|d$ or $D(r_0)|2^ld, 
D(r_0) \nmid 2^{l-1}d$ for some $1 \leq l \leq 
3$ and $j \geq l+\delta$.\\
We examine each case separately.
\vskip 5mm

\noindent
When $D(r_0)|d$, we have $\sqrt{r_0} \in \mathbb{Q}(\zeta_d)$ and so, 
$\sqrt{r_0} \in K'_j$. \\
Moreover, if $\delta\ge 1$, then
$[\mathbb{Q}(\zeta_{2^{{\rm min}(j,\delta)}}, \sqrt{r_0}): \mathbb{Q}]=
2^{{\rm min}(j,\delta)}=[K_j:\mathbb Q],$  except in the case when 
$\mathbb{Q}(\sqrt{r_0}) = \mathbb{Q}(\sqrt{2})$ 
which we have excluded in this section. Also, when $\delta=0$, 
$[\mathbb{Q}(\sqrt{r_0}): \mathbb{Q}] = 2 = [K'_j: \mathbb{Q}]$.
Therefore $K'_j = \mathbb{Q}(\zeta_{2^{{\rm min}(j,\delta)}}, \sqrt{r_0})$ if  
$D(r_0)|d.$

\noindent
When $D(r_0)|2^ld,  D(r_0) \nmid 2^{l-1}d$ for some $1 \leq l \leq 
3$ and $j \geq l+\delta$, it means that $D'|d'$, $\delta_0 = \delta+l$.
If $r_0 = u/v$, note that $\mathbb{Q}(\sqrt{r_0}) = 
\mathbb{Q}(\sqrt{uv}).$ 
Now, if $uv$ is odd, it has to be $\equiv 3({\rm mod~}4)$ since otherwise 
$D(r_0) = uv$ which cannot divide $2^ld$ without dividing $d$.
Also then $D(r_0) = 4uv = 4D'$, $D'|d'$, $\delta_0=2= \delta+l$ means 
that $l=1= \delta$ or $l=2, \delta =0.$
In case $uv \equiv 3({\rm mod~}4)$, we have $\sqrt{-r_0} \in 
\mathbb{Q}(\sqrt{d})$ as the discriminant of
$\mathbb{Q}(\sqrt{-r_0}) = -uv = D'$ which divides $d'$ and hence divides 
$d$.  Therefore 
$K'_j =  \mathbb{Q}(\zeta_{2^{{\rm min}(j,\delta)}}, \sqrt{-r_0})$,
when $D(r_0)|2^ld, D(r_0) \nmid 2^{l-1}d$ for some $1 \leq l 
\leq 3$ and $j \geq l+\delta$ and $r_0=u/v$ with $uv$ odd.
Here, we have used the fact that since $j \geq \delta_0 = 2$,  
$\zeta_4$ (and hence $\sqrt{-r_0}$) belongs to $L'_j$.
\vskip 4mm

\noindent
When $uv = 2s_0$ with $s_0 > 1$ odd, then $D(r_0) = 4uv=8s_0, \delta_0 = 
3, D'=s_0$. Also $\delta= \delta_0-l=3-l$ and $s_0=D'|d'$.
Thus, if $s_0 = \prod_{i=1}^k p_i$, then 
$\sqrt{t}\in\mathbb{Q}(\zeta_{p_1\cdots p_k}) \subseteq \mathbb{Q}(\zeta_d)$,
where $t:=\prod_{i=1}^k (\frac{-1}{p_i}) p_i$.
We have used the fact that $\sqrt{2}, i \in \mathbb{Q}(\zeta_8)$ and that 
$j \geq \delta_0 = 3$.
Hence 
$K'_j =  \mathbb{Q}(\zeta_{2^{{\rm min}(j,\delta)}}, \sqrt{t})$
when $uv$ is even and $D(r_0)|2^ld, 
D(r_0) \nmid 2^{l-1}d$ for some $1 \leq l \leq 3$ and $j \geq l+\delta$.
\qed
\vskip 5mm

\noindent
An immediate consequence of the previous lemma is the following result on 
the values of $\tau(j)$ and $\tau'(j)$.
\vskip 5mm

\begin{Lem}

If $j \leq \lambda+1$, then
$\tau(j)=1 \Leftrightarrow {\rm min}(j,\delta) \leq \gamma$.\\
If $j > \lambda+1$ and if either $D'\nmid d'$ or $\delta_0 > {\rm max}(j, 
\delta)$, then
$\tau(j)=1 \Leftrightarrow {\rm min}(j,\delta) \leq \gamma$.\\
If $j > \lambda + 1$ and $D(r_0)|d$, then 
$\tau(j)=1 \Leftrightarrow {\rm min}(j,\delta) \leq \gamma$ and 
$(\frac{D(r_0)}{c}) = 1$.\\
If $j > \lambda + 1$ and $D'|d', \delta < \delta_0 \leq j$ 
with $uv$ odd where $r_0 = u/v$, then
$\tau(j)=1 \Leftrightarrow {\rm min}(j,\delta) \leq \gamma$ and 
$(\frac{D(-r_0)}{c}) = 1$.\\
If $j > \lambda + 1$ and $D'|d', \delta < \delta_0 \leq j$ 
with $uv = 2 \prod_{i=1}^k p_i$ where $r_0 = u/v$ and $p_i$'s odd primes, 
then $\tau(j)=1 \Leftrightarrow {\rm min}(j,\delta) \leq \gamma$ and 
$(\frac{D(\prod_{i=1}^k (\frac{-1}{p_i})p_i)}{c}) = 1$.\\
We have $\tau'(j) = \tau(j)$ for $j \neq \lambda+1$.\\
If either $D' \nmid d'$ or $\delta_0 > {\rm max}(\lambda+1, 
\delta)$, then
$\tau'(\lambda+1)=1 \Leftrightarrow {\rm min}(\lambda+1,\delta) \leq \gamma$.\\
If $D(r_0)|d$, then 
$\tau'(\lambda+1)=1 \Leftrightarrow {\rm min}(\lambda+1,\delta) \leq \gamma$ and 
$(\frac{D(r_0)}{c}) = 1$.\\
If $D'|d', \delta < \delta_0 \leq \lambda+1$ 
with $uv$ odd where $r_0 = u/v$, then
$\tau'(\lambda+1)=1 \Leftrightarrow {\rm min}(\lambda+1,\delta) \leq \gamma$ and 
$(\frac{D(-r_0)}{c}) = 1$.\\
If $D'|d', \delta < \delta_0 \leq \lambda+1$ 
with $uv = 2 \prod_{i=1}^k p_i$
 where $r_0 = u/v$, then
$\tau'(\lambda+1)=1 \Leftrightarrow {\rm min}(\lambda+1,\delta) \leq 
\gamma$ and 
$( \frac{D(\prod_{i=1}^k (\frac{-1}{p_i})p_i)}{c} ) = 1$.
\end{Lem}

\section{Tables for the density $\delta_{a,b}(c,d)$ when
$\mathbb{Q}(\sqrt{r_0}) \neq \mathbb{Q}(\sqrt{2})$}
Recall that the density $\delta_{a,b}(c,d)$ is given by (\ref{begin}).
Since the primes considered are in $\phi(d)$ residue classes, it is more 
natural to compute the sum 
\begin{equation}
\label{begin2}
S := \phi(d) \delta_{a,b}(c,d) = \phi(d) \sum_{j \geq 1} \Big( 
\frac{\tau(j)}{[N_j:\mathbb{Q}]} -
\frac{\tau'(j)}{[N'_j:\mathbb{Q}]} \Big).
\end{equation}
Note that $S$ gives the
relative density of divisibility of $S_{a,b}$, that is
$$S=\lim_{x\rightarrow \infty}{\#\{p\le x:p\equiv c({\rm mod~}d),~p|S_{a,b}\}\over
\#\{p\le x:p\equiv c({\rm mod~}d)\}}.$$
Putting in the degrees of $N_j,N'_j$ we can simplify the sum in (\ref{begin2})  as follows.\\
\indent Since $[N_j:\mathbb{Q}] = [N'_j:\mathbb{Q}]$ and
$\tau(j) = \tau'(j)$ for $j \leq \lambda$, 
the terms corresponding to $j \leq \lambda$ do not contribute. Also
$\tau(j) = \tau'(j)$ for $j > \lambda$+1, but $\tau(\lambda+1)$ and 
$\tau'(\lambda+1)$ may be different (only) when 
$(C_{\lambda+1})$ holds. Therefore, we have :
\vskip 3mm

\begin{tabular}{llll}
$\frac{S}{\phi(2^{\delta})}$ & = 
$\tau(\lambda+1)2^{1-{\rm max}(\lambda+1,\delta)} -
\tau'(\lambda+1)2^{1-{\rm max}(\lambda+1,\delta)} $ & &\\
&+ $2^{\lambda+1}\sum_{j > \lambda+1, (C_j) {\rm ~fails}} 
\tau(j)2^{-{\rm max}(j,\delta)-j}$ &
& if $(C_{\lambda+1})$ holds\\
 &+ $
2^{\lambda+2}\sum_{j > \lambda+1, (C_j) {\rm ~holds}} 
\tau(j)2^{-{\rm max}(j,\delta)-j}$ & &

\end{tabular}
\vskip 3mm

\begin{tabular}{llll}
$\frac{S}{\phi(2^{\delta})}$ & = 
$\tau(\lambda+1)2^{1-{\rm max}(\lambda+1,\delta)} -
\tau(\lambda+1)2^{-{\rm max}(\lambda+1,\delta)} $ & &\\
&+ $2^{\lambda+1}\sum_{j > \lambda+1, (C_j) {\rm ~fails}} 
\tau(j)2^{-{\rm max}(j,\delta)-j}$ &
& if $(C_{\lambda+1})$ fails\\
 &+ $2^{\lambda+2}
\sum_{j > \lambda+1, (C_j) {\rm ~holds}} 
\tau(j)2^{-{\rm max}(j,\delta)-j}$ & &

\end{tabular}

\vskip 8mm

\noindent
As the degrees of the fields $N_j,N'_j$ and the values of 
$\tau(j),\tau'(j)$'s depend on the following three conditions, 
is convenient to have 3 tables depending on them. 
The three conditions are :\\ 
(A) $D'\nmid d'$;\\ 
(B) $D'|d', \delta_0 \leq 
\delta$;\\ 
(C) $D'|d', \delta_0 > \delta$.
\vskip 3mm

\noindent
Let us first work out the expression for $S$ in case A.\\
\vskip 2mm

\noindent
{\bf Case A: $D' \nmid d'$ }
\vskip 4mm

\noindent
Here, every $(C_j)$ fails. In particular,
$$\frac{S}{\phi(2^{\delta})}  = 
\tau(\lambda+1)2^{-{\rm max}(\lambda+1,\delta)}
+ 2^{\lambda+1}\sum_{j > \lambda+1} 
\tau(j)2^{-{\rm max}(j,\delta)-j}.$$
Moreover, since
$K_j=K'_j = \mathbb{Q}(\zeta_{2^{{\rm min}(j,\delta)}})$ for all $j \geq 
\lambda+1$, we have :\\
For all $j \geq \lambda+1$, $\tau(j) = \tau'(j)$ and
this is $1$ if and only if min$(j,\delta) \leq \gamma$.\\
Thus,
$S = \phi(2^{\delta}) 2^{\lambda+1}\sum_{j>\lambda, {\rm min}(j,\delta) \leq \gamma} 
2^{-{\rm max}(j,\delta)-j} = 
\phi(2^{\delta})2^{\lambda+1}(S_1+S_2)$,\\
where $S_1$ is the sum over $j \leq \delta$ and $S_2$ is the sum over $j 
\geq \delta+1$.\\
We get $$S_1 = \sum_{\lambda+1 \leq j\leq {\rm min}(\gamma,\delta)}2^{-\delta-j}{\rm ~and~}
S_2 = \cases{\sum_{j \geq {\rm max}(\lambda+1,\delta+1)}4^{-j} & if $\delta \leq \gamma$;\cr 
0 & otherwise.\cr}$$
From this, it is easy to obtain Table 1.
\vskip 5mm

\noindent
{\bf Case B : $D'|d'$, $\delta_0 \leq \delta$ }
\vskip 4mm

\noindent
Note that $(C_j)$ holds for all $j$.\\
Here $K_{\lambda+1} =  \mathbb{Q}(\zeta_{2^{{\rm min}(\lambda+1,\delta)}})$ 
and
$K'_{\lambda+1} = 
\mathbb{Q}(\zeta_{2^{{\rm min}(\lambda+1,\delta)}},\sqrt{r_0})$.\\ 
For all $j > \lambda+1$, we have
$K_j=K'_j  =  \mathbb{Q}(\zeta_{2^{{\rm min}(j,\delta)}}, \sqrt{r_0})$.\\
Therefore, $\tau(\lambda+1) = 1$ if and only if min$(\lambda+1,\delta)
\leq \gamma$;\\
$\tau'(\lambda+1) = 1$ if and only if min$(\lambda+1,\delta)
\leq \gamma$ and $(\frac{D(r_0)}{c}) = 1$.\\
Moreover, for $j > \lambda+1$, we have\\
$\tau(j)=\tau'(j)$ which is $1$ if and only if
min$(j,\delta) \leq \gamma$ 
and $(\frac{D(r_0)}{c}) = 1$.\\
Hence, we have

\begin{tabular}{llll}
$\frac{S}{\phi(2^{\delta})}$ & = 
$\tau(\lambda+1)2^{1-{\rm max}(\lambda+1,\delta)} -
\tau'(\lambda+1)2^{1-{\rm max}(\lambda+1,\delta)} $ & &\\
 &+ $2^{\lambda+2}\sum_{j > \lambda+1} 
\tau(j)2^{-{\rm max}(j,\delta)-j},$ & &
\end{tabular}\\
\noindent which  can be written down more explicitly as
$S = \phi(2^{\delta})(t_1+t_2+S_0)$, where\\
$$t_1 = \cases{2^{1-{\rm max}(\lambda+1,\delta)} &if  min$(\lambda+1,\delta) \leq \gamma$;\cr
0 & otherwise,\cr}$$ 
$$t_2 = \cases{- 2^{1-{\rm max}(\lambda+1,\delta)} & if 
min$(\lambda+1,\delta) \leq \gamma$ and $(\frac{D(r_0)}{c}) = 1$;\cr 
0 & otherwise,\cr}$$
$$S_0 = \cases{2^{\lambda+2}\sum_{j>\lambda+1,{\rm min}(j,\delta)\leq\gamma} 
2^{-{\rm max}(j,\delta)-j} & if
 $(\frac{D(r_0)}{c}) = 1$;\cr
0 & otherwise.\cr}$$
Further, $S_0 = S_{01} + S_{02}$, where $S_{01}$ is the subsum where $j$ 
varies over $j \leq \delta$ and $S_{02}$ is the subsum where $j$ varies 
over $j > \delta$. 
We find 
$$S_{01} = \cases{2^{\lambda+2-\delta}(2^{-1-\lambda}-2^{-{\rm 
min}(\gamma,\delta)}) & if $(\frac{D(r_0)}{c}) = 1$
and $\lambda+2 \leq {\rm min}(\delta,\gamma)$;\cr
0 & otherwise,\cr}$$
and that
$$S_{02} = \cases{2^{\lambda+2-2{\rm max}(\lambda+1,\delta)}/3
& if $(\frac{D(r_0)}{c}) = 1$ and $\delta \leq \gamma$;\cr
0 & otherwise.\cr}$$
From this, we obtain Table 2.
\vskip 5mm

\noindent
Finally, we work out the expression for $S$ in case C.
We  write $r_0=u/v$ and $t = -r_0$ or $\prod_{i=1}^k 
(\frac{-1}{p_i})p_i$ according as to whether $uv$ is odd or 
$uv = 2 \prod_{i=1}^k p_i$. We also write $D(t)$ for the 
discriminant of the quadratic field ${\mathbb Q}(\sqrt{t})$.
\vskip 5mm

\noindent
{\bf Case C : $D'|d'$, $\delta_0 > \delta$}
\vskip 4mm

\noindent
Notice that there are finitely many $j$'s for which the property $(C_j)$ 
may fail in this case. 
Now $$K_{\lambda+1} =  
\mathbb{Q}(\zeta_{2^{{\rm min}(\lambda+1,\delta)}}),~K'_{\lambda+1} = 
\cases{\mathbb{Q}(\zeta_{2^{{\rm min}(\lambda+1,\delta)}}) & if 
$\lambda+1 < \delta_0$; \cr
\mathbb{Q}(\zeta_{2^{{\rm min}(\lambda+1,\delta)}},\sqrt{t}) & otherwise.}$$
For all $j>\lambda+1$, we have
$$K_j=K'_j  =\cases{ \mathbb{Q}(\zeta_{2^{{\rm min}(j,\delta)}}) & if  
$j < \delta_0$;\cr
\mathbb{Q}(\zeta_{2^{{\rm min}(j,\delta)}}, \sqrt{t}) & otherwise.\cr}$$
So, we have $\tau(\lambda+1) = 1$ if and only if min$(\lambda+1,\delta)
\leq \gamma$ and furthermore we have
$$\tau'(\lambda+1) = 1 \iff 
\cases{{\rm min}(\lambda+1,\delta)
\leq \gamma,~\lambda+1 < \delta_0 ; & \cr 
{\rm min}(\lambda+1,\delta) 
\leq \gamma,~\lambda+1 \geq \delta_0,{\rm ~and~} 
(\frac{D(t)}{c}) = 1. & \cr}$$
Moreover, for $j > \lambda+1$ with $j < \delta_0$, we 
have $\tau(j)=\tau'(j)$ which is $1$ if and only if
min$(j,\delta) \leq \gamma$.On the other hand, for $j > \lambda+1$ with 
$j \geq \delta_0$, 
we have $\tau(j)=\tau'(j)$ which is $1$ if and only if min$(j,\delta) 
\leq \gamma$ 
and $(\frac{D(t)}{c}) = 1$.\\
Therefore, we get $S = \phi(2^{\delta})(t_1+t_2+ S_1+S_2)$, where
$$t_1 = \tau(\lambda+1)2^{1-{\rm max}(\lambda+1,\delta)};$$
$$t_2 =\cases{ -\tau'(\lambda+1)2^{1-{\rm max}(\lambda+1,\delta)} & if  
 $\lambda+1 \geq \delta_0$;\cr
-\tau(\lambda+1)2^{-{\rm max}(\lambda+1,\delta)} & if
 $\lambda+1 < \delta_0$;\cr}$$
$S_1 = 2^{\lambda+1}
\sum \{ 2^{-{\rm max}(j,\delta)-j} : j>\lambda+1, 
~j \delta_0,{\rm min}(j,\delta)\leq \gamma \}$;\\
$S_2 = 2^{\lambda+2}\sum \{2^{-{\rm max}(j,\delta)-j} : j>\lambda+1, 
j \geq \delta_0,{\rm min}(j,\delta)\leq \gamma \}$ if
$(\frac{D(t)}{c})=1$ and, is $0$, otherwise.\\
\indent Putting in the values of $\tau(\lambda+1)$ and $\tau'(\lambda+1)$, we 
obtain
$$ t_1 =\cases{2^{1-{\rm max}(\lambda+1,\delta)} & if 
min$(\lambda+1,\delta) \leq \gamma$;\cr 
0 & otherwise,\cr}$$
$$t_2 =\cases{ -2^{1-{\rm max}(\lambda+1,\delta)} & if  
$\lambda+1 \geq \delta_0,{\rm 
min}(\lambda+1,\delta)\leq \gamma, 
(\frac{D(t)}{c})=1$;\cr
-2^{-{\rm max}(\lambda+1,\delta)} & if
 $\lambda+1 < \delta_0,{\rm min}(\lambda+1,\delta)\leq \gamma$;\cr
0 & otherwise.\cr}$$
Finally, as before, we break up each of $S_1$ and $S_2$ into two subsums 
over $j \leq \delta$, respectively, over $j > \delta$.
So, we have $S_1 = S_{11} + S_{12}$, where
$$S_{11} = 2^{\lambda+1-\delta} \sum \{2^{-j} : {\rm min}(\gamma,\delta) 
\geq j > \lambda+1 \};$$
$$S_{12} = \cases{2^{\lambda+1} \sum \{4^{-j} : \delta_0>j \geq {\rm 
max}(\lambda+2, \delta+1) \} & if $\delta \leq \gamma$;\cr
0 & otherwise.\cr}$$
Similarly, we have $S_2 = S_{21}+S_{22}$,where
$$S_{21} =0,~S_{22} = \cases{2^{\lambda+2} \sum \{4^{-j} : j \geq {\rm 
max}(\lambda+2,\delta_0) \} & if $\delta \leq \gamma$ and 
$(\frac{D(t)}{c})=1$;\cr 
0 & otherwise.}$$ 
On evaluating these expressions further we obtain Table 3.

\section{The intersection fields when $\mathbb{Q}(\sqrt{r_0}) = 
\mathbb{Q}(\sqrt{2})$}
Next we consider the case where $r_0 = 2$ or $1/2$. 
Note that the discriminant of $\mathbb{Q}(\sqrt{2})$ is $8$ and that 
$\sqrt{2}$ belongs to the cyclotomic field
$\mathbb{Q}(\zeta_8)$ (indeed $\sqrt{2} = \zeta_8 + \zeta_8^{-1}$). Also
note that $\mathbb Q(i,\sqrt{2})=\mathbb Q(\zeta_8)$ (we have
$\zeta_8=(i+1)/\sqrt{2}$).
For $j \geq 1$ we consider as before the degrees of the fields 
$N_j,N'_j$.
The earlier expressions in Lemma \ref{tweegraad} are valid and, in fact, simplify to 
give:
\begin{Lem}
The degrees of $N_j = \mathbb{Q}(\zeta_{2^j},r^{1/2^{j-1}},\zeta_d)$
and $N'_j = \mathbb{Q}(\zeta_{2^j},r^{1/2^{j}},\zeta_d)$ over $\mathbb{Q}$ 
are given by :
$$\frac{1}{\phi(d')}[N_j: \mathbb{Q}] = 
\cases{
2^{{\rm max}(j,\delta)-1} & if $j\leq \lambda + 1$;\cr
 2^{{\rm max}(j,\delta)+j-\lambda-3} &  if
$j > \lambda+1$ and $3 \leq {\rm max}(j,\delta)$;\cr
2^{{\rm max}(j,\delta)+j-\lambda-2} &  if
$j> \lambda+1$ and $3 > {\rm max}(j,\delta)$,}
$$
$$\frac{1}{\phi(d')}[N'_j: \mathbb{Q}]=\cases{ 
2^{{\rm max}(j,\delta)-1}  & if $j \leq \lambda$;\cr
2^{{\rm max}(j,\delta)+j-\lambda-2} & if
$j> \lambda$ and $3 \leq {\rm max}(j,\delta)$;\cr
2^{{\rm max}(j,\delta)+j-\lambda-1}  & if
$j> \lambda$ and $3>{\rm max}(j,\delta)$.}$$
\end{Lem}

\noindent
The fields 
$K_j =  \mathbb{Q}(\zeta_{2^j},r^{1/2^{j-1}}) \cap \mathbb{Q}(\zeta_d)$
and
$K'_j = \mathbb{Q}(\zeta_{2^j},r^{1/2^j}) \cap \mathbb{Q}(\zeta_d)$
are to be determined. 
This is where the computation gives different values from 
Lemma \ref{vvier}. However, the method of evaluation is the same and
the degrees turn out to be :
\vskip 5mm

\noindent
For $j > \lambda+1$, 
$$[K_j : \mathbb{Q}] = \cases{2^2 & if $j \leq 2, \delta \geq 3$;\cr
2^{{\rm min}(j,\delta)-1} & if either $j \geq 3,~\delta\ge 1$ or $j<3,~1 \leq \delta \leq 2$;\cr 
1 & if $\delta = 0$.}$$

\noindent
For $j> \lambda$,
$$
[K'_j : \mathbb{Q}] =\cases{2^j & if $j \leq 2, \delta \geq 3$;\cr
2^{{\rm min}(j,\delta)-1} & if either $j \geq 3,~\delta\ge 1$ or $j<3,~1 \leq \delta \leq 2$;\cr
1 & if $\delta = 0$.}
$$

\noindent
As we have evidently, $K_j = 
\mathbb{Q}(\zeta_{2^{{\rm min}(j,\delta)}})$ 
for $j \leq \lambda+1$ and 
for every $j$, $\mathbb{Q}(\zeta_{2^{{\rm min}(j,\delta)}})$ is a subfield of 
$K_j$, we have  the following result :
\vskip 3mm

\begin{Lem} We have
$K_j = \mathbb{Q}(\zeta_{2^{{\rm min}(j,\delta)}})$ for all $j$ unless 
$\lambda=0,j=2,\delta\geq3$.\\
In the exceptional cases $\lambda=0,j=2,\delta\geq3$,
we have $K_2 = \mathbb{Q}(\zeta_{2^{{\rm min}(j,\delta)}},
\sqrt{2})=\mathbb{Q}(i,\sqrt{2}) = \mathbb{Q}(\zeta_8)$.\\
\indent Further, we have $K'_j = \mathbb{Q}(\zeta_{2^{{\rm min}(j,\delta)}})$ for 
all $j$ unless 
$\lambda<j\leq 2,\delta\geq3$. The exceptional cases here are : either $\lambda=0,j=1,\delta \geq 3$ or 
$\lambda \leq 1,j=2,\delta \geq 3$. 
We find the following intersection fields:
$$\cases{\lambda=0,j=1,\delta\geq 3, & 
$K'_1 = \mathbb{Q}(\zeta_{2^{{\rm min}(j,\delta)}},\sqrt{2})=\mathbb{Q}(\sqrt{2})$;\cr
\lambda \leq 1,j=2,\delta\geq 3, & 
$K'_2 = \mathbb{Q}(\zeta_{2^{{\rm min}(j,\delta)}},\sqrt{2})
=\mathbb{Q}(i,\sqrt{2}) =  \mathbb{Q}(\zeta_8)$.\cr}$$
\end{Lem}

\section{Tables for the density when
$\mathbb{Q}(\sqrt{r_0}) = \mathbb{Q}(\sqrt{2})$}

\noindent
Let $S$ be defined as in (\ref{begin2}). We divide its computation into four cases :\\ 
(A) $\delta < 3$;\\
(B) $\delta \geq 3$ and $\lambda \geq 2$,\\ 
(C) $\delta \geq 3$ and $\lambda =1$, and\\
(D) $\delta \geq 3$ and $\lambda =0$.
\vskip 5mm

\noindent
{\bf Case A : $\delta < 3$}
\vskip 5mm

\noindent
Then $K_j = K'_j = \mathbb{Q}(\zeta_{2^{{\rm min}(j,\delta)}})$ for all $j$.
Thus $\tau(j)=\tau'(j)$ for all $j$ and, this is $1$ if and only if
min$(j,\delta) \leq \gamma$. It turns out that
$S = \phi(2^{\delta})(t_1+t_2+t_3)$, with
$$t_1 =\cases{2^{-{\rm max}(\lambda+1,\delta)} & if $\lambda \leq 1, 
{\rm min}(\lambda+1,\delta) \leq \gamma$;\cr 
0 & otherwise,\cr},~t_2 =\cases{1/8 & if $\lambda=0, \delta \leq \gamma$;\cr 
0 & otherwise,}$$
$$t_3 = \cases{2^{\lambda+2-2{\rm max}(\lambda+1,2)}/3 & if $\delta \leq \gamma$;\cr 
0 & otherwise,}$$
where $t_1,t_2,t_3$ correspond, respectively, to the terms in (\ref{begin2})
with $j=\lambda+1$, $\lambda+2\le j\le 3$, $j\ge {\rm max}(3,\lambda+2)$ 
and $j\ge {\rm max}(3,\delta+1)$. From this, we obtain Table 4.\\

\noindent
{\bf Case B : $\delta \geq 3$, $\lambda \geq 2$}
\vskip 5mm

\noindent
Once again,
$K_j = K'_j = \mathbb{Q}(\zeta_{2^{{\rm min}(j,\delta)}})$ for all $j$. 
Note that $(C_j)$ always holds true.
We obtain $$S=\varphi(d)\sum_{j\ge \lambda+2\atop {\rm min}(j,\delta)\le \gamma}
\Big({1\over [N_j:\mathbb Q]}-{1\over [N_j':\mathbb Q]}\Big)=
\phi(2^{\delta})(t_1+t_2),$$ where 
$$t_1 =\cases{ 
2^{1-\delta}-2^{\lambda+2-\delta-{\rm min}(\gamma,\delta)} & if $\lambda+2 \leq 
{\rm min}(\gamma,\delta)$;\cr 
0 & otherwise,\cr}$$
$$t_2 =\cases{2^{\lambda+2-2{\rm max}(\lambda+1,\delta)}/3 & if $\delta \leq 
\gamma$;\cr 
0 & otherwise,\cr}$$
with $\varphi(2^{\delta})t_1$, $\varphi(2^{\delta})t_2$ the subsum over
$j\le \delta$, respectively $j>\delta$.\\

\noindent
{\bf Case C : $\delta \geq 3$, $\lambda =1$}
\vskip 5mm

\noindent
Here, we need to observe that 
when $8|d$, the Galois automorphism
$\zeta_d \mapsto \zeta_d^c$ of $\mathbb{Q}(\zeta_d)$ fixes
$\sqrt{2}$ if and only if $c \equiv \pm 1({\rm mod~}8)$.
We obtain 
$${S\over \varphi(2^{\delta})}={\tau(\lambda+1)\over 2^{\delta-1}}-
{\tau'(\lambda+1)\over 2^{\delta-1}}
+2^{\lambda+2}\sum_{3\le j\le \delta}{\tau(j)\over 2^{{\rm max}(j,\delta)+j}}+
2^{\lambda+2}\sum_{j>{\rm max}(2,\delta)}
{\tau(j)\over 2^{{\rm max}(j,\delta)+j}},$$
which can be written as $t_1+t_2+t_3+t_4$ say, where further evaluation yields that
$$t_1 = \cases{2^{1-\delta} & if $2 \leq \gamma$;\cr
0 & otherwise;\cr},~t_2 =\cases{-2^{1-\delta} & if $3 \leq \gamma$;\cr
0 & otherwise;\cr}$$
$$t_3 =\cases{ 2^{1-\delta} - 2^{3-\delta-{\rm min}(\gamma,\delta)} & if $3 \leq \gamma$;\cr
0 & otherwise;\cr},~{\rm and~}t_4 =\cases{2^{3-2\delta}/3 & if $\delta \leq \gamma$;\cr 
0 & otherwise.\cr}$$
Table 5 is obtained from cases B and C.
\vskip 5mm

\noindent
{\bf Case D : $\delta \geq 3$, $\lambda =0$}
\vskip 5mm

\noindent
As in the previous case, we need the fact that
when $8|d$, the Galois automorphism
$\zeta_d \mapsto \zeta_d^c$ of $\mathbb{Q}(\zeta_d)$ fixes
$\sqrt{2}$ if and only if $c \equiv \pm 1 ({\rm mod~ }8)$.\\
We find that $S = \phi(2^{\delta})(t_1+t_2+t_3+t_4)$, where
$$t_1 =\cases{2^{1-\delta} & if $c\equiv \pm 3({\rm mod~}8)$;\cr
0 & otherwise,\cr},~t_2 =\cases{2^{-\delta} & if $3 \leq \gamma$; \cr
0 & otherwise,\cr}$$
$$t_3 = \cases{2^{-\delta}-2^{2-\delta-{\rm min}(\gamma,\delta)} & if $3 \leq 
{\rm min}(\gamma,\delta)$;\cr 
0 & otherwise,\cr},~{\rm and~}t_4 =\cases{2^{2-2\delta}/3 & if $\delta \leq \gamma$;\cr
0 & otherwise,}$$
where $t_1,t_2,t_3,t_4$ correspond, respectively, to the terms in (\ref{begin2})
with $j=1$, $j=2$, $3\le j\le \delta$ and $j\ge {\rm max}(3,\delta+1)$. 
This yields us Table 6.

\section{Extremal densities}
We have $0\le \varphi(d)\delta_{a,b}(c,d)\le 1$. In this section we are 
interested
when $\delta_{a,b}(c,d)=0$ and when $\delta_{a,b}(c,d)=1/\varphi(d)$. The 
following
elementary result shows that if $c\not\equiv 1({\rm 
mod~}(d,2^{\lambda+1}))$, then $\delta_{a,b}(c,d)=0$.
\begin{Lem}
If $p\nmid (a,b)$ and $p|S_{a,b}$, then $p\equiv 1({\rm mod~}2^{\lambda+1})$.
\end{Lem}
{\it Proof}. For a prime $p$ put $\tau(p)=(p-1)/(p-1,h)$.
If $p\nmid (a,b)$ and $p|ab$, then $p\nmid S_{a,b}$, so we may assume
that $p\nmid ab$. Since $r^{\tau(p)}=(r_0^h)^{\tau(p)}\equiv 1({\rm mod~}p)$ by Fermat's
little theorem, it follows that ord$_p(r)|\tau(p)$. 
If $p$ is to divide
$S_{a,b}$, then $\tau(p)$ must be even and so $\nu_2(p-1)\ge \lambda + 1$. \qed
\begin{Thm} 
\label{moetnaamhebben}
{\rm a)} Suppose that $\delta_{a,b}(c,d)=0$.
This happens if and only if\\
i) $\lambda \ge \gamma$ and $\delta > \gamma$;\\
or\\
ii) $\lambda=\gamma-1$, $\delta > \gamma$, $D(r_0)|d$ and 
$({D(r_0)\over c})=1$.\\
Moreover, if $\delta_{a,b}(c,d)=0$, then there are
at most finitely primes $p\equiv c({\rm mod~}d)$ dividing the
sequence $S_{a,b}$.\\
{\rm b)} Suppose that $\delta_{a,b}(c,d)=1/\varphi(d)$. This happens if and only if\\
i) $\lambda=0$, $\delta=0$, $D(r_0)|d$ and $({D(r_0)\over c})=-1$;\\
or\\
ii) min$(\gamma,\delta) > \lambda$, $D(r_0)|d$ and 
$({D(r_0)\over 
c})=-1$.\\
Moreover, if $\delta_{a,b}(c,d)=1/\varphi(d)$, then there are
at most finitely primes $p\equiv c({\rm mod~}d)$ not dividing the
sequence $S_{a,b}$.
\end{Thm}
{\it Proof}. 
For a prime $p$ put $\tau(p)=(p-1)/(p-1,h)$. The first parts 
of both (a) and (b) follow on 
inspection of the Tables. Let us prove the second part of (a) now.
If $\lambda \geq \gamma$ and $\delta > \gamma$, we claim that $\tau(p)$ 
 is odd. Indeed, 
writing $p = c + qd$, and $c-1 = 
2^{\gamma}c_0$ with $c_0$ odd, we have
$p-1 = 2^{\gamma}c_0 + 2^{\delta}qd'$.
Therefore, $v_2(p-1) = \gamma$ since $\delta > \gamma$.
Now, $(p-1,h) = (p-1,2^{\lambda}h')$ which has $2$-adic valuation $\gamma$ 
since $\lambda \geq \gamma$. Therefore $\tau(p)$ is odd in the case (i) 
of (a) of the theorem. Since clearly
ord$_p(r)|\tau(p)$, it then follows that $p\nmid S_{a,b}$.
Finally suppose we are in case ii. Suppose that $p>2$ is a prime satisfying
$p\equiv c({\rm mod~}d)$ and such that $p$ does not divide $ab$. Then, by 
the properties of the Kronecker symbol,
$$({{\overline{r_0}}\over p})=\Big({D(r_0)\over p}\Big)=\Big({D(r_0)\over c}\Big)=1,$$
where the first symbol is the Legendre symbol and ${\overline{r_0}}$ denotes the
reduction of $r_0$ modulo $p$.
It follows that
$$r_0^{h(p-1)\over 2(p-1,h)}\equiv 1({\rm mod~}p),$$
and so ord$_p(r)|\tau(p)/2$.
We claim that $\tau(p)/2$ is odd. Now
$p-1 = 2^{\gamma}c_0 + 2^{\delta}qd'$ which has $2$-adic valuation $\gamma$ because $\delta > \gamma$.
On the other hand, $2(p-1,h) = 2(p-1,2^{\lambda}h')= 
2(p-1,2^{\gamma-1}h')$ which has $2$-adic valuation $1 + (\gamma 
-1) = \gamma$. Thus, $\tau(p)/2$ is odd and so $p\nmid S_{a,b}$. \\
b) The proof is similar; let us consider (i) first.\\
As $\delta = \lambda = 0$, we have $h$ is odd and $r = r_0^h$.
If $p>2$ is a prime not dividing $ab$, then
$$({\overline{r_0}\over p})=\Big({D(r_0)\over p}\Big)=\Big({D(r_0)\over c}\Big)=-1$$
by assumption. Thus, $r_0^{(p-1)/2} \equiv -1({\rm mod`}p)$, which implies that 
$r^{(p-1)/2} \equiv -1({\rm mod~}p)$ and therefore, that $p|S_{a,b}$. Finally
suppose we are in case ii.
Writing $p = c + qd$, and $c-1 = 
2^{\gamma}c_0$ with $c_0$ odd, we have
$p-1 = 2^{\gamma}c_0 + 2^{\delta}qd'$.
Therefore, $v_2(p-1) \geq {\rm min}(\delta,\gamma)$.
Now, $v_2(p-1,h) = v_2(p-1,2^{\lambda}h') = \lambda$, since  
$v_2(p-1) \geq {\rm min}(\gamma,\delta) > \lambda$. 
Therefore, we have that $\frac{h}{(p-1,h)}$ is odd while 
$\tau(p)$ is even; that is,
${\frac{p-1}{2(p-1,h)}}$ is a positive integer.  
Once again, we have for each prime not dividing $2ab$ that
$$({\overline{r_0}\over p})=\Big({D(r_0)\over p}\Big)=\Big({D(r_0)\over c}\Big)=-1.$$
Thus,  $(r_0^{(p-1)/2})^{\frac{h}{(p-1,h)}} \equiv -1({\rm mod~}p)$. 
But then $r^{\frac{p-1}{2(p-1,h)}} = (r_0^{(p-1)/2})^{\frac{h}{(p-1,h)}} 
\equiv -1 ({\rm mod~}p)$, which means that $p|S_{a,b}$. \qed\\

\noindent {\tt Example}. 1) By case ii of (a) we infer that $\delta_{3,1}(11,12)=0$ 
(cf. Conjecture 1.1 of Fermat).\\
2) By case ii of (b) we infer that  $\varphi(8)\delta_{2,1}(\pm 3,8)=1$ (easily
proved using $(2/p)=(-1)^{(p^2-1)/8}$), cf. the paper by 
Sierpi\'nski \cite{sier}.\\

\noindent Perhaps a more illuminating phrasing of the above theorem is the 
following.
\begin{Thm}
For a prime $p$ put $\tau(p)=(p-1)/(p-1,h)$.\\
{\rm a)}  We have $\delta_{a,b}(c,d)=0$ if and only if $\tau(p)$ is odd or
$2||\tau(p)$ and $({r_0\over p})=1$, for all but finitely many primes 
$p\equiv c({\rm mod~}d)$.\\
{\rm b)} We have $\delta_{a,b}(c,d)=1/\varphi(d)$ if and only if for all but finitely
many primes $p\equiv c({\rm mod~}d)$ we have that $\tau(p)$ is even and
$({r_0\over p})=-1$.
\end{Thm}

\noindent {\tt Conclusion}: if the density is extremal, then this can always be explained by
elementary arguments not using more than quadratic reciprocity and, furthermore,
the associated set of exceptional primes is at most finite.\\

\noindent {\tt Remark 5} ({\it uniform distribution}). It is generally not true 
that the primes dividing $S_{a,b}$ are 
uniformly distributed over the residue classes modulo $d$. However, there are
some cases where we have uniform distribution. For example, if $d$ is odd and $D(r_0)\nmid d$, then 
the primes in any residue class mod $d$ which divide 
$S_{a,b}$ have the same density.

\section{Some numerical experiments}
For each entry in Tables 1-6 an example
with parameters $a$ and $b=1$ was choosen and below we give the
value of $\delta_{a,1}(c,d)$ according to the tables
on the one hand, and an approximation to this that
consists of the first six decimals of the ratio
$${\#\{p\le p_m: p\equiv c({\rm mod~}d),~p|S_{a,1}\}\over \#\{p\le p_m:p\equiv c({\rm
mod~}d)\}},$$
where $p_m$ denotes the $m$th prime and $m=2097152000\approx 2\cdot 10^9$.
As a rule of thumb an approximation of $\delta_{a,1}(c,d)$
obtained in this way by looking for prime divisors amongs
the primes should have an accuracy of about $\pi(p_m;d,c)^{-1/2}$. We
clearly observed in our experiments that for larger $d$ the accuracy
tends to be less (and the same holds for  the run time).\\

\noindent{\bf ~~~~Test cases for Table 1}
\vskip 5mm

\begin{tabular}{|l|l|l|l|}
\hline
Residue class & $a$ & $\phi(d) \delta_{a,1}(c,d)$ 
& Experimental value\\ \hline 
\hline
$17$ mod $56$ & $3^2$ & $5/6$ &$0.833200\cdots$\\
\hline
$17$ mod $56$ & $3^8$ & $1/3$ & $0.333317\cdots$\\
\hline
$1~$ mod $21$ & $5$ & $2/3$ & $0.666592\cdots$\\
\hline
$7~$ mod $20$ & $3^4$  & $0$ & 0\\
\hline
$7~$ mod $20$ & $3^3$ & $1/2$ & $0.500015\cdots$\\
\hline
\end{tabular}
\vskip 5mm

\noindent{\bf ~~~~Test cases for Table 2}
\vskip 5mm

\begin{tabular}{|l|l|l|l|} \hline
Residue class & $a$ & $\phi(d) \delta_{a,1}(c,d)$ 
 & Experimental value\\ \hline 
\hline
$9~$ mod $28$ & $7^2$ & $1/3$ & $0.333312\cdots$\\
\hline
$5~$ mod $12$ & $3^2$ & $1$ & 1\\
\hline
$1~$ mod $15$ & $5$ & $1/3$ & $0.333257\cdots$\\
\hline
$7~$ mod $15$ & $5$ & $1$ & 1\\
\hline
$1~$ mod $12$ & $3$ & $2/3$ & $0.666657\cdots$\\
\hline
$5~$ mod $12$ & $3$ & $1$ & 1\\
\hline
$11$ mod $20$ & $5^4$ & $0$ & 0\\
\hline
$13$ mod $24$ & $3$ & $1/2$ & $0.500006\cdots$\\
\hline
$13$ mod $56$ & $7$ & $1$ & 1\\
\hline
$7~$ mod $20$ & $5^2$ & $0$ & 0\\
\hline
\end{tabular}
\vskip 5mm

\noindent {\bf ~~~~Test cases for Table 3}
\vskip 5mm

\begin{tabular}{|l|l|l|l|} \hline
Residue class  & $a$ & $\phi(d) \delta_{a,1}(c,d)$ 
 & Experimental value\\ \hline 
\hline
$1~$ mod $12$ & $6$ & $11/12$ & $0.916693\cdots$\\
\hline
$5~$ mod $12$ & $6$ & $3/4$ & $0.749989\cdots$\\
\hline
$1~$ mod $12$ & $6^2$ & $5/6$ & $0.833362\cdots$\\
\hline
$5~$ mod $12$ & $6^2$ & $1/2$ & $0.499996\cdots$\\
\hline
$7~$ mod $12$ & $6$ & $1/2$ & $0.500038\cdots$\\
\hline
$11$ mod $28$ & $14^2$ & $0$ & 0\\
\hline
$7~$ mod $12$ & $6^4$ & $0$ & 0\\
\hline
$7~$ mod $30$ & $6^2$ & $5/12$ & $0.416679\cdots$\\
\hline
$11$ mod $30$ & $6^2$ & $1/4$ & $0.250055\cdots$\\
\hline
$7~$ mod $30$ & $6^4$ & $1/12$ & $0.083321\cdots$\\
\hline
$11$ mod $30$ & $6^4$ & $1/4$ & $0.250055\cdots$\\
\hline
$7~$ mod $15$ & $6$ & $17/24$ & $0.708336\cdots$\\
\hline
$11$ mod $15$ & $6$ & $5/8$ & $0.624999\cdots$\\
\hline
$7~$ mod $15$ & $6^4$ & $1/12$ & $0.083321\cdots$\\
\hline
$11$ mod $15$ & $6^4$ & $1/4$ & $0.250055\cdots$\\
\hline
\end{tabular}

\noindent {\bf ~~~~Test cases for Table 4}
\vskip 5mm

\begin{tabular}{|l|l|l|l|} \hline
Residue class & $a$ & $\phi(d) \delta_{a,1}(c,d)$ 
& Experimental value\\ \hline 
\hline
$5$ mod $14$ & $2$ & $17/24$ & $0.708327\cdots$\\
\hline
$5$ mod $12$ & $2$ & $11/12$ & $0.916652\cdots$\\
\hline
$7$ mod $12$ & $2$ & $1/2$ & $0.499961\cdots$\\
\hline
$7$ mod $12$ & $2^2$ & $0$ & 0\\
\hline
$5$ mod $6$ & $2^2$ & $5/12$ & $0.416673\cdots$\\
\hline
$5$ mod $12$ & $2^2$ & $5/6$ & $0.833331\cdots$\\
\hline
$5$ mod $6$ & $2^8$ & $1/24$ & $0.041672\cdots$\\
\hline 
$5$ mod $12$ & $2^4$ & $1/6$ & $0.166685\cdots$\\
\hline 
$7$ mod $12$ & $2^4$ & $0$ & 0\\
\hline 
\end{tabular}
\vskip 5mm

\noindent {\bf ~~~~Test cases for Table 5}
\vskip 5mm

\begin{tabular}{|l|l|l|l|} \hline
Residue class  & $a$ & $\phi(d) \delta_{a,1}(c,d)$ 
 & Experimental value\\ \hline 
\hline
$5~$ mod $24$ & $2^4$ & $0$ & 0\\
\hline
$17$ mod $24$ & $2^4$ & $1/3$ & $0.333372\cdots$\\
\hline
$17$ mod $48$ & $2^4$ & $2/3$ & $0.666740\cdots$\\
\hline
$17$ mod $96$ & $2^4$ & $1/2$ & $0.500145\cdots$\\
\hline
$41$ mod $48$ & $2^4$ & $0$ & 0\\
\hline
$17$ mod $24$ & $2^2$ & $2/3$ & $0.666659\cdots$\\
\hline
$7$ mod $24$ & $2^2$ & $0$ & 0\\
\hline
$5~$ mod $24$ & $2^2$ & $1$ & 1\\
\hline
$17$ mod $32$ & $2^2$ & $3/4$ & $0.750049\cdots$\\
\hline
\end{tabular}
\vskip 5mm

\noindent {\bf ~~~~Test cases for Table 6}
\vskip 5mm

\begin{tabular}{|l|l|l|l|} \hline
Residue class & $a$ & $\phi(d) \delta_{a,1}(c,d)$ 
 & Experimental value\\ \hline 
\hline
$9$ mod $40$ & $2$ & $5/6$ & $0.833411\cdots$\\
\hline
$7$ mod $8$ & $2$ & $0$ & 0\\
\hline
$5$ mod $8$ & $2$ & $1$ & 1\\
\hline
$9$ mod $16$ & $2$ & $3/4$ & $0.749983\cdots$\\
\hline
\end{tabular}\\
\vskip 5mm

\noindent {\tt Acknowledgement}. This paper was written during a stay
February-March 2007 of the second author at the Max-Planck-Institut f\"{u}r
Mathematik. The authors have the pleasure in thanking that institute 
for providing excellent hospitality and a wonderful work
atmosphere. In addition they thank Yves Gallot for kindly writing a 
Visual $C^{++}$ program
that was used to create the data in the test case tables.

\vfil\eject

\medskip\noindent {\footnotesize Max-Planck-Institut f\"ur Mathematik,\\
Vivatsgasse 7, D-53111 Bonn, Germany.\\
e-mail: {\tt moree@mpim-bonn.mpg.de}}
\vskip 5mm

\medskip\noindent {\footnotesize Statistics \& Mathematics Unit, 
Indian Statistical Institute,\\ 8th Mile Mysore Road, Bangalore 560059, India.\\
email : {\tt sury@isibang.ac.in}}

\end{document}